\documentclass[10pt]{amsart}
\usepackage[dvips,ps,color,all]{xy}
\usepackage{graphicx,amssymb,latexsym,amsfonts,verbatim,amscd,amsmath,amsthm}
\usepackage[colorlinks=true,linkcolor=black]{hyperref}
\providecommand{\dis}{\displaystyle}

\newtheorem{thm}{Theorem}[section]
\newtheorem{lemma}[thm]{Lemma}
\newtheorem{rem}[thm]{Remark}
\newtheorem{rems}[thm]{Remarks}

\newtheorem{cor}[thm]{Corollary}

\newtheorem{prop}[thm]{Proposition}
\newtheorem{ex}[thm]{Example}
\newtheorem{exs}[thm]{Examples}

\begin{document}
\title{On Borel fixed ideals generated in one degree}
\author{Achilleas Sinefakopoulos}
\begin{abstract}
We construct a (shellable) polyhedral cell complex that supports a
minimal free resolution of a Borel fixed ideal, which is minimally
generated (in the Borel sense) by just one monomial in
$S=\Bbbk[x_1,x_2,...,x_n]$; this includes the case of powers of
the homogeneous maximal ideal $(x_1,x_2,...,x_n)$ as a special
case.
   In our most general result we prove that for any Borel fixed ideal
$I$ generated in one degree, there exists a polyhedral cell
complex that supports a minimal free resolution of $I$.
\end{abstract}
\date{February 20, 2007}
\maketitle

    \section{Introduction} \label{S:I}

We study resolutions over the polynomial ring
$S=\Bbbk[x_1,x_2,...,x_n]$, where $\Bbbk $ is a field. The idea to
encode the structure of the resolution of a monomial ideal in the
combinatorial structure of
 a  simplicial complex was introduced  in
\cite{DPS} (see also \cite{E2}). The idea was generalized later in
\cite{DS}, where resolutions supported on a
\textbf{\emph{regular}} cell complex were introduced. The
generalization continued in \cite{Batz} and \cite{Wel}, where
monomial resolutions supported on a CW-complex were introduced and
studied. Recently,  the necessity for CW-resolutions is justified
in \cite{Velasco}, and their sufficiency is disproved by the
existence of a monomial ideal whose resolution cannot be supported
on a CW-complex.

In this paper we study Borel fixed ideals generated in the same
degree $d$, which we call $d$-generated. For $d$-generated Borel
fixed ideals a minimal free resolution is already well known,
namely the Eliahou-Kervaire resolution (see, e.g. \cite{Elia} or
\cite{Miller}), which is also CW-cellular, as it is proved in
\cite{Batz} by using discrete Morse theory. Moreover, in
\cite{Batz}, the authors give the Morse complex that supports a
minimal free resolution for powers of the homogeneous maximal
ideal $(x_1,...,x_n)$ of the polynomial ring
$S=\Bbbk[x_1,x_2,...,x_n]$ as a worked example. More generally,
the Morse complex that supports a minimal free resolution of
principal Borel fixed ideals, that is, of those ideals which are
minimally generated (in the Borel sense) by just one monomial, is
given in \cite{Wel}. However, it is not clear whether any of those
Morse complexes is regular or not. Thus, a natural question is
whether there exists a \textbf{\emph{regular}} cell complex that
supports a minimal resolution of a $d$-generated Borel fixed
ideal. We answer the question positively in this paper, which is
organized as follows:

 In Section 2, we give the basic notation and preliminaries
   for the rest of this paper and we refer to the literature for more details.

  In Section 3, we answer the above  natural question by constructing inductively a
   shellable polyhedral cell complex that supports the minimal free
   resolution of a principal Borel fixed ideal in $S=\Bbbk[x_1,x_2,...,x_n]$;
this includes the case of powers of the homogeneous maximal ideal
as a special case. Our most general result is theorem \ref{P3:CR},
where we prove that for any
   $d$-generated Borel fixed ideal $I$ , there exists a polyhedral cell
complex  that supports a minimal free resolution of $I$. it should
be noted that  the basis we use
   in the minimal free resolution is different than the one used in the Elliahou-Kervaire resolution.

 Finally, in  Section 4, we consider the lcm-lattice of
   a  $d$-generated Borel fixed ideal. In particular,
   in proposition \ref{K:LCM}, we show that it is
    ranked. This result was proved (in greater generality) independently in \cite{Phan}.

    \section{Notation-Preliminaries}\label{S:NP}
        \subsection{Monomial ideals}

All ideals in this paper are considered to be \emph{monomial}
ideals. We work over the polynomial ring $S=\Bbbk[x_1,x_2,...,x_n]$
with $char(\Bbbk)=0$. For small $n$
    we may use the letters $a,b,c,d,...$ instead of $x_1,x_2,x_3,x_4,...$.

For a monomial $\mathbf{m}=x_1^{a_1}x_2^{a_2}...x_n^{a_n}$ in $S$,
we define the exponent vector to be
    $e(\mathbf{m})=(a_1,a_2,...,a_n)$ and we set
    $\text{max}(\mathbf{m})$ to be the largest index of a variable that
    divides $\mathbf{m}$.

We let $G(I)$ denote the unique minimal set of monomial generators
of a (monomial) ideal $I$. A (monomial) ideal $I$ is called
\emph{Borel fixed}, if for every $\mathbf{m}$ in $G(I)$ and every
$x_t$ that divides $\mathbf{m}$,
\[
          \mathbf{m}_{t \to s}:=\frac{\mathbf{m}}{x_{t}}x_{s}
            \]
is in $I$ for all $1\leq s <t$. A Borel fixed ideal $I$ is called
\textbf{principal} Borel, and it is written as
\[
   I=<\mathbf{m}>,
\]
if $I$ is the smallest Borel fixed ideal such that $\mathbf{m}$ is
in $G(I)$. In this case, we also say that $I$ is generated by just
one monomial $\mathbf{m}$ in the Borel sense.

\begin{ex} \label{EX:NP} Let $S=\Bbbk[a,b,c]$. The ideal $(a^2,a b, b^2, ac,bc)$ is
a Borel fixed ideal, which is also principal Borel, because
\[
        (a^2,a b, b^2, ac,bc)=<bc>
  \]
\end{ex}

For more on monomial ideals we refer to \cite{E1}, \cite{E2} and
\cite{Miller}.

\subsection{Cellular resolutions and polyhedral complexes}

    As in \cite{DS}, let $X$ be a \emph{regular} cell complex having
$G(I)$, the set of minimal generators of $I$, as its set of vertices
and let $\epsilon_X$ be an \emph{incidence function} on $X$. It is
well known that such a function exists, (see e.g. pp. 244-248 in
\cite{Massey}). Next we label each nonempty face $F$ of $X$ by the
least common multiple $\mathbf{m}_F$ of the monomials $\mathbf{m}_j$
in $G(I)$, which correspond to the vertices of $F$. The
\emph{degree} $\mathbf{a}_F$ of the face $F$ is defined to be the
exponent vector $e(\mathbf{m}_F)$.

Let $SF$ be the free $S$-module with one generator $F$ in degree
$\mathbf{a}_F$. The \emph{cellular complex} $\mathbf{F}_X$ is
     the $\mathbb{Z}^n$-graded $S$-module $\underset{\varnothing \ne F \in
     X}{\bigoplus} SF$ with differential
\[
  \partial F= \underset{\varnothing \ne F' \in X}{\sum} \epsilon_X
  (F,F ') \frac{\mathbf{m}_F}{\mathbf{m}_{F'}} F'.
\]
For each degree $\mathbf{b} \in \mathbb{Z}^n$ let $X_{\preceq
\mathbf{b}}$ be the subcomplex of $X$ on the vertices of degree
$\preceq \mathbf{b}$. The following results are proved in \cite{DS},

\begin{prop}\label{P1:NP} The complex $\mathbf{F}_{X}$ is a free resolution of
$I$ if and only if $X_{\preceq \mathbf{b}}$ is acyclic over $\Bbbk$
for all degrees $\mathbf{b}$. In this case, $\mathbf{F}_{X}$ is
called a cellular resolution of $I$.
\end{prop}

\begin{cor} \label{C1:NP}
 The cellular complex $\mathbf{F}_{X}$ is a resolution of $I$ if and
 only if the cellular complex $\mathbf{F}_{X_{\preceq \mathbf{b}}}$
 is a resolution of the monomial ideal $I_{\preceq \mathbf{b}}$ for
 all $\mathbf{b} \in \mathbb{Z}^n$.
\end{cor}

\begin{rem}\label{R1:NP} A cellular resolution $\mathbf{F}_{X}$ is minimal if
and only if any two comparable faces $F' \subseteq F$ of the same
degree coincide.
\end{rem}

\begin{ex} Let $I \subset S$ be a monomial ideal with $G(I)=\{x_1^d,x_2^d,...,x_n^d\}$
for a fixed positive integer $d$. Then the labelled
$(n-1)$-simplex $\Delta_{n-1} (x_1^d,x_2^d,...,x_n^d)$ with
vertices in $G(I)$ supports a minimal free resolution of $I$.
\end{ex}

Note that in this paper, whenever we say \textbf{cellular}
resolutions, we mean resolutions supported on a \textbf{regular}
cell complex. Otherwise, we talk about CW-resolutions to emphasize
the difference and avoid confusion.

 The above results are presented in \cite{Miller} for
polyhedral complexes, which is a special case of regular cell
complexes. A \emph{polyhedral cell complex} $X$ is a nonempty
finite collection of convex polytopes (in some real vector space
$\mathbb{R}^{N}$), called faces of $X$, satisfying two properties:
 \begin{itemize}
    \item If $P$ is a polytope in $X$ and $F$ is a face of $P$, then
    $F$ is in $X$.
    \item If $P$ and $Q$ are in $X$, then $P \cap Q$ is a face of
    both $P$ and $Q$.
 \end{itemize}

 For more on polytopal complexes we refer to \cite{Ziegler}. Here it suffices to
mention that a basic notion that we are going to use is that of
\emph{regular subdivisions} of a polytope (\cite{Ziegler}, p.129,
or \cite{Bruns}, p.34). Another concept is the \emph{shellability}
of a polytope (\cite{Ziegler}, p.233).

\subsection{Results from Algebraic Topology}

We assume familiarity with the basic notions of $CW$-complex and
\textbf{regular} cell complex and their differences. Recall that
the closures of the cells of a regular CW-complex are homeomorphic
with closed balls. For example, any polyhedral cell complex is
regular. So we only state the two major theorems from algebraic
topology that we use. We need the cellular version of
Mayer-Vietoris theorem and the K\"unneth theorem with field
coefficients.

 \begin{thm}[Mayer-Vietoris]\label{MV:NP}. Let $X$ be a $CW$-complex
 and let $Y_1$ and $Y_2$ be $CW$ subcomplexes of $X$ such that $X=Y_1 \cup Y_2$. Then
 there is an exact sequence
\[
      \cdots \rightarrow {\tilde H}_i (Y_1 \cap Y_2 ; \Bbbk) \rightarrow {\tilde H}_i (Y_1 ; \Bbbk) \oplus {\tilde H}_i (Y_2 ; \Bbbk)\rightarrow {\tilde H}_i (X ; \Bbbk)\rightarrow{\tilde H}_{i-1} (Y_1 \cap Y_2 ; \Bbbk) \rightarrow
      \cdots.
                            \]
\end{thm}

\begin{thm} [K\"unneth]\label{K:NP} Let $X$ and $Y$ be two $CW$-complexes. Then
there is a natural isomorphism
\[
   \bigoplus_{j} \left(H_j (X ; \Bbbk) \otimes_{\Bbbk}H_{i-j} (Y;\Bbbk)\right) \rightarrow H_i(X \times Y
   ;\Bbbk).
   \]
\end{thm}
We refer to \cite{Hatcher} or \cite{Massey},  for more on these.

    \section{Cellular Resolutions of $d$-generated Borel fixed ideals}\label{S:CR}
\subsection{Three basic Lemmas}
Now we may proceed to our study of $d$-generated Borel fixed
ideals. Let $I$ and $J$ be two monomial ideals in $S$ and assume
that $X$ and $Y$ are regular cell complexes in $\mathbb{R}^{N}$
(for some $N$) that support a (minimal) free resolution of  $I$
and $J$, respectively.

Can we say anything about the cellular resolution of $I+J$ and/or
the cellular resolution of $IJ$?

The following three lemmas give some results related to this
question, which will be useful in proving our main results. The
assumption that the cell complexes are regular is not necessary.

\begin{lemma}\label{L1:CR} Let $I$ and $J$ be two monomial ideals in
$S$ such that $G(I+J)=G(I) \cup G(J)$ set-theoretically. Suppose
that
\begin{itemize}
 \item[(i)] $X$ and $Y$ are regular cell complexes in some
 $\mathbb{R}^{N}$ that support a (minimal) free resolution of
            $I$ and $J$, respectively, and
 \item[(ii)] $X\cap Y$ is a regular cell complex that supports a
            (minimal) free resolution of $I \cap J$.
            \end{itemize}
Then $X\cup Y$ supports a (minimal) free resolution of $I+J$.
\end{lemma}

\textbf{Proof}: First let $Z:=X \cup Y$ and note that $Z$ is a
regular cell complex. From our hypothesis, we have
\[
  {\tilde H}_i (X_{\preceq \mathbf{b}} ; \Bbbk)=0, \quad {\tilde H}_i (Y_{\preceq \mathbf{b}} ;
  \Bbbk)=0, \quad \text{and} \quad {\tilde H}_i ((X \cap Y)_{\preceq \mathbf{b}} ; \Bbbk)=0
\]
for all $i$ and all $\mathbf{b} \in \mathbb{Z}^n$. Furthermore, it
is clear that
\[
  Z_{\preceq \mathbf{b}}=(X \cup Y)_{\preceq \mathbf{b}}= X_{\preceq \mathbf{b}}\cup Y_{\preceq \mathbf{b}}
\]
for all $\mathbf{b} \in \mathbb{Z}^n$. Then the Mayer-Vietoris
 theorem \ref{MV:NP} gives us the following exact sequence
\[
      {\tilde H}_i (X_{\preceq \mathbf{b}} ; \Bbbk) \oplus {\tilde H}_i (Y_{\preceq \mathbf{b}} ; \Bbbk)\rightarrow {\tilde H}_i (Z_{\preceq \mathbf{b}} ; \Bbbk)\rightarrow{\tilde H}_{i-1} ((X \cap Y)_{\preceq \mathbf{b}} ; \Bbbk)
\]
Consequently, ${\tilde H}_i (Z_{\preceq \mathbf{b}}; \Bbbk)=0$ and
the proof is complete from proposition \ref{P1:NP}. \endproof

\begin{rem} For any two monomial
ideals $I$ and $J$, we have
\[
    G(I+J) \subseteq G(I) \cup G(J).
\]
Our assumption that $G(I+J)=G(I) \cup G(J)$ guarantees the right
labelling of the cell complex $X\cup Y$. A case where equality
becomes true is when all elements of $G(I)\cup G(J)$ are of the
same degree.

Note that from the labelling of $X$,$Y$ and $X\cap Y$ and our
assumptions above, it follows that
\[
    G(I \cap J)= G(I) \cap G(J).
    \]

\end{rem}

\begin{lemma}\label{L2:CR} Let $I \subset \Bbbk [x_1,...,x_k]$ and $J
\subset \Bbbk [x_{k+1},...,x_n]$ be two monomial ideals. Suppose
that  $X$ and $Y$ are regular cell complexes  in some
$\mathbb{R}^{N}$ of dimension $k-1$ and $n-k-1$, respectively,
that support a (minimal) free resolution of $I$ and $J$,
respectively. Then the regular cell complex $X \times Y$ supports
a (minimal) free resolution for $IJ$.
\end{lemma}
\textbf{Proof}: Let $Z:=X \times Y$ and let
$\mathbf{b}=(\mathbf{b}_1,\mathbf{b}_2) \in \mathbb{Z}^n$, where
$\mathbf{b}_1 \in \mathbb{Z}^k$ and $\mathbf{b}_2 \in
\mathbb{Z}^{n-k}$. Then, it is easy to check that
\[
  Z_{\preceq \mathbf{b}}=(X \times Y)_{\preceq \mathbf{b}}= X_{\preceq \mathbf{b}_1}\times Y_{\preceq \mathbf{b}_2}
\]
From the K\"unneth theorem \ref{K:NP} for $CW$ complexes, there is
an isomorphism
\[
   \bigoplus_{j} \left(H_j (X_{\preceq \mathbf{b}_1} ; \Bbbk) \otimes_{\Bbbk}H_{i-j} (Y_{\preceq
   \mathbf{b}_2};\Bbbk)\right) \cong H_i(X_{\preceq \mathbf{b}_1}\times Y_{\preceq \mathbf{b}_2}  ;
   \Bbbk)= H_i (Z_{\preceq \mathbf{b}} ; \Bbbk)
\]
for all $i$. From our hypothesis, we have
\[
  {\tilde H}_i (X_{\preceq \mathbf{b}_1} ; \Bbbk)=0 \quad \text{and} \quad {\tilde H}_i ( Y_{\preceq \mathbf{b}_2} ;
  \Bbbk)=0,
\]
for all $i$. Therefore,
\[
  H_0 (Z_{\preceq \mathbf{b}} ; \Bbbk) =\Bbbk
  \otimes_{\Bbbk} \Bbbk=\Bbbk,
\]
while
\[
 H_i (Z_{\preceq \mathbf{b}} ;\Bbbk) =0,
\]
for $i>0$. Now assume that $e_{X} \times e_{Y}$ and $\sigma_{X}
\times \sigma_{Y}$ are two comparable faces of $X\times Y$ with the
same label. That is,
\[
  e_{X} \subset \sigma_{X} \quad \text{and} \quad  e_{Y} \subset \sigma_{Y}
\]
and
\[
  label(e_{X} \times e_{Y})=label(\sigma_{X}\times \sigma_{Y})=(\mathbf{b}_1,\mathbf{b}_2)
\]
where $\mathbf{b}_1 \in \mathbb{Z}^k$ and $\mathbf{b}_2 \in
\mathbb{Z}^{n-k}$. Then,
\[
 label(e_{X})=label(\sigma_{X})=(\mathbf{b}_1)
\]
and
\[
  label(e_{Y})=label(\sigma_{Y})=(\mathbf{b}_2).
\]
Therefore, $e_X=\sigma_X$ and $e_Y=\sigma_Y$, $e_{X} \times
e_{Y}=\sigma_{X} \times \sigma_{Y}$. The proof is complete from
proposition \ref{P1:NP} and remark \ref{R1:NP}.
\endproof

\begin{rem}\label{remark:CR} From our conclusion in lemma \ref{L2:CR}, it follows that
   \begin{align}\notag
   pdim(S/IJ)&=dim(X \times Y)+1\\\notag
             &=(k-1)+(n-k-1) +1\\\notag
             &=n-1.\notag
    \end{align}
\end{rem}

\begin{lemma}\label{L3:CR}  Let $I \subset \Bbbk [x_1,...,x_k]$ and $J
\subset \Bbbk [x_{k},...,x_n]$ be two monomial ideals such that
$|G(IJ)|=|G(I)|\cdot |G(J)|$. Suppose that there exists a regular
cell complex $X$ in $\mathbb{R}^{N}$ (for some $N$) of dimension
$k-1$ and a regular cell complex $Y$ in $\mathbb{R}^{n-k}$ of
dimension $n-k$ , which support a (minimal) free resolution of $I$
and $J$ respectively. Then the regular cell complex $X \times Y$
supports a (minimal) free resolution of $IJ$.
\end{lemma}
 \textbf{Proof}: Let $Z:=X \times Y$ and let
$\mathbf{b}=(\mathbf{b}_1,\beta,\mathbf{b}_2) \in \mathbb{Z}^n$,
where $\mathbf{b}_1 \in \mathbb{Z}^{k-1}$, $\beta \in \mathbb{Z}$
and $\mathbf{b}_2 \in \mathbb{Z}^{n-k}$. Then, for $1\leq k \leq
\beta -1$ define iteratively $Z_{\preceq \mathbf{b}}^{(k)}$ as
follows
\[
 Z_{\preceq \mathbf{b}}^{(k+1)}=Z_{\preceq \mathbf{b}}^{(k)}\cup \left(X_{\preceq (\mathbf{b}_1,k)}\times Y_{\preceq
  (\beta -k,\mathbf{b}_2)}\right).
\]
where $Z_{\preceq \mathbf{b}}^{(1)}=X_{\preceq
(\mathbf{b}_1,0)}\times Y_{\preceq (\beta,\mathbf{b}_2)}$, and note
that
 \[
 Z_{\preceq \mathbf{b}}=Z_{\preceq \mathbf{b}}^{(\beta)}=\left( X_{\preceq (\mathbf{b}_1,0)}\times Y_{\preceq
  (\beta,\mathbf{b}_2)}\right) \cup \left(X_{\preceq (\mathbf{b}_1,1)}\times Y_{\preceq
  (\beta -1,\mathbf{b}_2)}\right) \cup \dots \cup \left(X_{\preceq (\mathbf{b}_1,k)}\times Y_{\preceq
  (\beta -k,\mathbf{b}_2)}\right)
\]
Moreover,
\[
 \left( X_{\preceq (\mathbf{b}_1,0)}\times Y_{\preceq
  (\beta,\mathbf{b}_2)}\right) \cap \left(X_{\preceq (\mathbf{b}_1,1)}\times Y_{\preceq
  (\beta -1,\mathbf{b}_2)}\right)= X_{\preceq (\mathbf{b}_1,0)}\times Y_{\preceq
  (\beta -1,\mathbf{b}_2)},
\]
or more generally,
\[
 Z_{\preceq \mathbf{b}}^{(k)} \cap \left(X_{\preceq (\mathbf{b}_1,k+1)}\times Y_{\preceq
  (\beta -k-1,\mathbf{b}_2)}\right)=X_{\preceq (\mathbf{b}_1,0)}\times Y_{\preceq
  (\beta -k-1,\mathbf{b}_2)}.
\]
 Therefore, by combining the Mayer-Vietoris theorem \ref{MV:NP} with the
 K\"unneth formula \ref{K:NP} we get
 \[
 {\tilde H}_i (Z_{\preceq \mathbf{b}}; \Bbbk)=0.
 \]
 Now assume that $e_{X} \times e_{Y}$ and $\sigma_{X}
\times \sigma_{Y}$ are two comparable faces of $X\times Y$ with the
same label. That is,
\[
  e_{X} \subset \sigma_{X} \quad \text{and} \quad  e_{Y} \subset \sigma_{Y}
\]
and
\[
  label(e_{X} \times e_{Y})=label(\sigma_{X}\times \sigma_{Y})=(\mathbf{b}_1,\beta,\mathbf{b}_2)
\]
where $\mathbf{b}_1 \in \mathbb{Z}^{k-1}$, $\beta \in \mathbb{Z}$
and $\mathbf{b}_2 \in \mathbb{Z}^{n-k}$. Then,
\[
 label(e_{X})=(\mathbf{b}_1,\beta_1) \quad \text{and} \quad
 label(\sigma_{X})=(\mathbf{b}_1,\beta_2),
\]
which implies $\beta_1 \leq \beta_2$ and
\[
  label(e_{Y})=(\beta-\beta_1,\mathbf{b}_2) \quad \text{and} \quad label(\sigma_{Y})=(\beta- \beta_2,\mathbf{b}_2),
\]
which implies $\beta-\beta_1 \leq \beta-\beta_2$, that is, $\beta_2
\leq \beta_1$. Thus we have $\beta_1 = \beta_2$ and then,
\[
 label(e_{X})=label(\sigma_{X})=(\mathbf{b}_1,\beta_1)
\]
and
\[
  label(e_{Y})=label(\sigma_{Y})=(\beta-\beta_1,\mathbf{b}_2).
\]
Therefore, $e_X=\sigma_X$ and $e_Y=\sigma_Y$, and so $e_{X} \times
e_{Y}=\sigma_{X} \times \sigma_{Y}$. Thus, the resolution is minimal
and the proof is complete.\endproof

\begin{rems} 1) For any two monomial
ideals $I$ and $J$, we have
\[
        G(IJ) \subseteq G(I)G(J).
\]
Thus, our assumption that $|G(IJ)|=|G(I)|\cdot |G(J)|$ forces
$G(IJ)=G(I)G(J)$.

2) Let $\mathbf{F_{X}}$ be the cellular resolution of $I$ and let
$\mathbf{F_{Y}}$ be the cellular resolution of $J$. Then define
\[
  \mathbf{F_{X\times Y}}:=\mathbf{F_{X}}\otimes \mathbf{F_{Y}}.
\]
(see e.g. pp. 280-282 in \cite {Massey}).
 \begin{ex}
    Let $S=\Bbbk[a,b,c]$. The resolution of $I=(a,b)$ is of the form
\[
 0\rightarrow S(-2) \rightarrow S^2(-1) \rightarrow (a,b)\rightarrow 0
                       \]
and the resolution of $J=(b,c)$ is of the form

\[
 0\rightarrow S(-2) \rightarrow S^2(-1) \rightarrow (b,c)\rightarrow 0
                       \]
Therefore, the resolution of $IJ=(a,b)(b,c)$ is of the form
\[
  0\rightarrow S(-4) \rightarrow S^4(-3) \rightarrow S^4(-2)\rightarrow IJ \rightarrow 0,
                       \]
which is the tensor product of the first two resolutions.
 \end{ex}

3) As in remark \ref{remark:CR}, from our conclusion in lemma 
\ref{L3:CR}, it follows that
   \begin{align}\notag
   pdim(S/IJ)&=dim(X \times Y)+1\\\notag
             &=(k-1)+(n-k) +1\\\notag
             &=n.\notag
    \end{align}

4) A lemma similar to lemmas \ref{L2:CR} and \ref{L3:CR} for
monomial ideals
 \[
    I \subset \Bbbk [x_1,...,x_{k-1},x_k] \quad \text{and} \quad  J \subset \Bbbk [x_{k-1},x_k,...,x_n]
    \]
    and corresponding regular cell complexes $X$ and $Y$ with
\[
    dim(X)=k-1 \quad \text{and} \quad dim(Y)=n-k+1
    \]
    would fail because we would have
    \begin{align}\notag
    dim(X \times Y)+1&=(k-1)+(n-k+1)+1\\\notag
                     &=n+1\\\notag
                     & >pdim(S/IJ).\notag
    \end{align}
    \end{rems}

\subsection{Powers of the homogeneous maximal ideal}

Now we may prove our first main result, which is about the powers
of the homogeneous maximal ideal in $S$.

\begin{thm} \label{P1:CR} There exists a (shellable)
polyhedral cell complex $P_d(x_1,...,x_n)$ that supports a minimal
free resolution of $(x_1,...,x_n)^d$. Moreover, $P_d(x_1,...,x_n)$
is a polyhedral subdivision of the $(n-1)$-simplex $\Delta_{n-1}
(x_1^d,x_2^d,...,x_n^d)$.
\end{thm}

\textbf{Proof}: The proof will be by induction on $d$. It is clear
that if $d=1$, then the standard $(n-1)$-simplex denoted by
$\Delta_{n-1} (x_1,x_2,...,x_n)$,  supports a minimal free
resolution of $(x_1,...,x_n)$ for all $n\geq 1$. Thus
\[
 P_1(x_1,...,x_n)= \Delta_{n-1}(x_1,x_2,...,x_n)
\]
for all $n\geq 1$. Also, $P_1(x_{k+1},...,x_n)$ is a subcomplex of
$P_1(x_{k},...,x_n)$ for all $k<n$. Next, assume that for some
$d\geq 1$ we have constructed $P_d(x_1,...,x_n)$ for all $n\geq 1$
and that $P_d(x_{k+1},...,x_n)$ is a subcomplex of
$P_d(x_{k},...,x_n)$ for all $k<n$. . Define the ideals
\[
   I_k= (x_1,x_2,...,x_k)(x_k,x_{k+1},...,x_n)^d
\]
and note that an easy (finite) induction on $k$ gives us
   \[
  I_1 +\cdots + I_{k}=(x_1,...,x_k)(x_1,x_2,...,x_n)^{d}.\\
 \]
 for all $1\leq k \leq n$. Indeed, assuming that we have proved it
 for $k-1$, for some $k>1$, then we have
\begin{align}\notag
    I_1 +\cdots+I_{k-1}+ I_{k}&=(x_1,...,x_{k-1})(x_1,x_2,...,x_n)^{d}+(x_1,...,x_k)(x_k,...,x_n)^d\\\notag
                                 &=(x_1,...,x_{k-1})(x_1,x_2,...,x_n)^{d}+x_k(x_k,...,x_n)^d\\\notag
                                 &=(x_1,...,x_{k-1})(x_1,x_2,...,x_n)^{d}+x_k(x_1,...,x_n)^d\\\notag
                                 &=(x_1,...,x_k)(x_1,...,x_n)^{d}\notag
\end{align}
Moreover, we see that
   \begin{align}\notag
 (I_1 + \cdots + I_{k})\cap I_{k+1}&= (x_1,...,x_k)(x_1,...,x_n)^{d}\cap (x_1,...,x_{k+1})(x_{k+1},...,x_n)^d\\\notag
           &=(x_1,...,x_k)(x_{k+1},...,x_n)^d.\notag
\end{align}
       From lemmas \ref{L2:CR} and \ref{L3:CR}, we conclude that the polyhedral cell
complexes $C_k$ and $D_k$ ($k=1,2,...,n$) defined by
\[
   C_k:=\Delta_{k-1}(x_1,x_2,...,x_k) \times P_d(x_k,...,x_n)
\]
and
\[
 D_k:=C_k \cap C_{k+1}=\Delta_{k-1}(x_1,x_2,...,x_k) \times P_d(x_{k+1},...,x_n)
\]
 support a  minimal free resolution for $I_k$ and $(I_1+\cdots+I_k) \cap I_{k+1}$,
 respectively.

 Thus, from this and lemma \ref{L1:CR}, the polyhedral cell complex
$C'_k$, which is defined recursively by
 \[
  C'_1=C_1,\quad \text{and} \quad C'_{k+1}=C'_{k} \cup C_{k+1}
 \]
 for $k\geq 1$, supports a (minimal) free resolution for $(x_1,...,x_k)(x_1,x_2,...,x_n)^{d}$.
 Accordingly, set
\[
  P_{d+1} (x_1,x_2,...,x_n):=C'_n= C_1 \cup C_2 \cup \cdots \cup C_n .
\]
and the construction of our polyhedral cell complex is done by
induction. The fact that $P_d(x_1,...,x_n)$ is a polyhedral
subdivision of the $(n-1)$-simplex $\Delta_{n-1}
(x_1^d,x_2^d,...,x_n^d)$ is clear from our construction.
$P_d(x_1,...,x_n)$ is a regular subdivision of the $(n-1)$-simplex
$\Delta_{n-1} (x_1^d,x_2^d,...,x_n^d)$ (see, e.g. \cite{Bruns},
p.37). Since a regular subdivision of a polytope is shellable
(\cite{Ziegler}, p.243), we conclude that $P_d(x_1,...,x_n)$ is
shellable.
\endproof

\begin{ex}
 Let $I=(a,b,c,d)^2$.  Using the software package MACAULAY 2 \cite{Mike}, we see that the
 polyhedral cell complex that supports the minimal free resolution
 of $I$ is
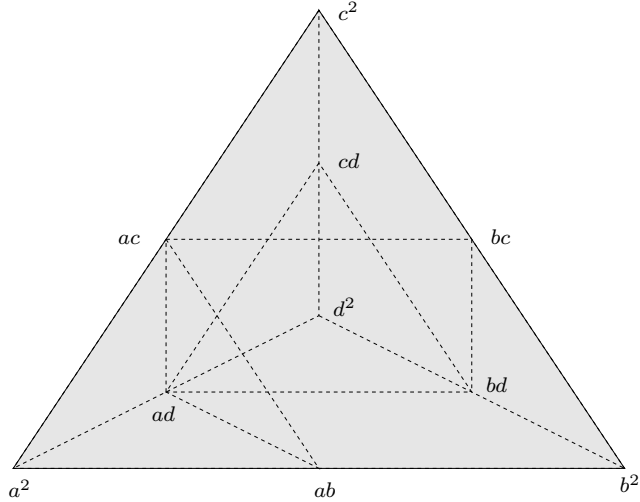
\begin{figure}[htbp]
\begin{center}
\input{maximal1.pstex_t}
\end{center}
\end{figure}

This can be decomposed as follows:

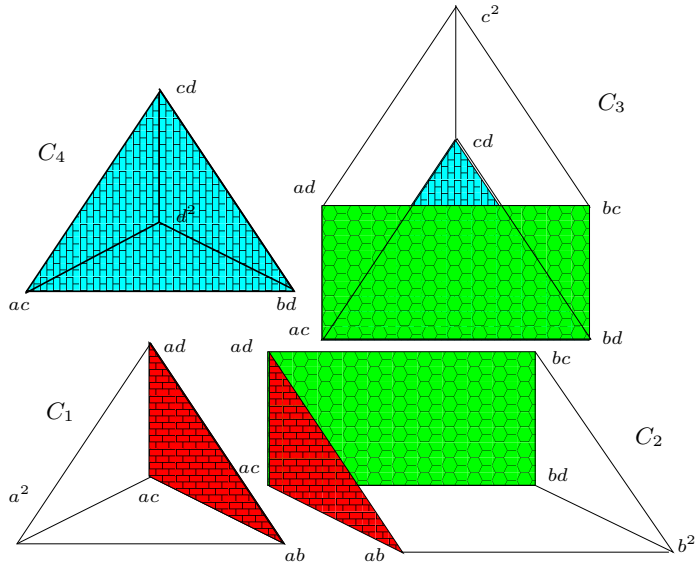
\begin{figure}[htbp]
\begin{center}
\input{maximal1parts.pstex_t}
\end{center}
\end{figure}
\end{ex}
Another cell complex that supports a minimal free resolution of
$I$ is the following (Morse complex), which supports the
Eliahou-Kervaire resolution of $I$. In particular, note that it is
not polyhedral.

\begin{figure}[htbp]
\begin{center}

\input{maximal1EK.pstex_t}

\end{center}
\end{figure}

\begin{rem} From theorem \ref{P1:CR} and corollary \ref{C1:NP}, we may get
 a minimal \emph{cellular} resolution for all ideals of the form $I_{\preceq \mathbf{b}}$
 ($\mathbf{b} \in \mathbb{Z}^n$) (see also \cite{PV}).
\end{rem}

\subsection{Principal Borel fixed ideals}
Our next goal is to prove a more general result for principal
Borel fixed ideals. Note that the following theorem includes
theorem \ref{P1:CR} as a special case, since
\[
  (x_1,...,x_n)^d=<x_n^d>.
\]

\begin{thm}\label{P2:CR} There exists a (shellable) polyhedral cell
complex $Q(\mathbf{m})$ that supports a minimal free resolution of
the principal Borel fixed ideal
\[
    I=<\mathbf{m}>= \prod_{j=1}^{s} I_{\lambda_j}^{d_j},
\]
where $\mathbf{m}=x_{\lambda_1}^{d_1}x_{\lambda_2}^{d_s}\dots
x_{\lambda_s}^{d_s}$, $I_i=(x_1,x_2,...,x_i)$ and $1\leq
\lambda_1<\lambda_2<...<\lambda_s$. Moreover, $Q(\mathbf{m})$ is a
subcomplex of $P_d(x_1,...,x_n)$, where $d=degree(\mathbf{m})$. In
particular, $Q(\mathbf{m})$ is the union of all the convex polytopes
(i.e. the faces) of the polyhedral cell complex $P_d(x_1,...,x_n)$,
with vertices in $<\mathbf{m}>$.
\end{thm}

\begin{rem}
 If $s=1$, then
$\mathbf{m}=x_{\lambda_{1}}^{d_1}$, and so
$Q(\mathbf{m})=P_{d_1}(x_1,x_2,...,x_{\lambda_{1}})$. If
$\lambda_{s-1}=1$, then $s=2$ and $\mathbf{m}=x_1^{d_1}
x_{\lambda_{2}}^{d_2}$, so $Q(\mathbf{m})$ is obtained by
multiplying all the labels of the vertices of
$P_{d_2}(x_1,x_2,...,x_{\lambda_{2}})$ by $x_1^{d_1}$.
\end{rem}

Before we prove the above theorem we need a lemma. Because of the
above remark, we may assume that $s>1$ and $\lambda_{s-1}>1$.

\begin{lemma}\label{L4:CR} Let $I$ be a principal Borel fixed ideal as
above. Define the ideals
\[
  N_k=<x_{\lambda_1}^{d_1}x_{\lambda_2}^{d_2}\cdots
  x_{\lambda_j}^{d_j}x_{k}^{d_{j+1}+...+d_{s-1}}>
\]
for $\lambda_j<k\leq \lambda_{j+1}$ ($j<s-1$ and $\lambda_0=0$).
Then for $\lambda_{s-1} < \mu\leq \lambda_s$
    \begin{itemize}
    \item[(a)]$
  N_i(x_1,...,x_\mu)^{d_s}=N_1
  (x_1,...,x_\mu)^{d_s}+N_2(x_2,...,x_\mu)^{d_s}+...+N_{i}(x_{i},...,x_\mu)^{d_s}$
\end{itemize}
  for $i=1,2...,\lambda_{s-1}$, and
  \begin{itemize}
      \item[(b)] $  N_j(x_1,...,x_\mu)^{d_s} \cap
  N_{j+1}(x_{j+1},...,x_\mu)^{d_s}=N_{j}(x_{j+1},...,x_\mu)^{d_s}$
  \end{itemize}
for $j=1,2...,\lambda_{s-1}-1$,

\end{lemma}

 \textbf{Proof}: (a) First it is clear that
\[
  N_i(x_1,...,x_\mu)^{d_s} \supseteq N_1
  (x_1,...,x_\mu)^{d_s}+N_2(x_2,...,x_\mu)^{d_s}+...+N_{i}(x_{i},...,x_\mu)^{d_s}.
  \]
  Next, note that the ideal
 $N_i(x_1,...,x_\mu)^{d_s}$ is principal Borel, minimally
 generated (in the Borel sense) by
 \[
    x_{\lambda_1}^{d_1}x_{\lambda_2}^{d_2}\cdots x_{\lambda_j}^{d_j}x_{i}^{d_{j+1}+...+d_{s-1}}
x_{\mu}^{d_s},
  \]
($\lambda_j<i\leq \lambda_{j+1}<\mu$, since $j<s-1$), which is
contained in $N_i(x_i,...,x_\mu)^{d_s}$. Thus, to complete the
proof of this part, it suffices to show that the sum of the ideals
on the right hand side of the equality to be proved is Borel
fixed. Set $J_k=N_k(x_k,...,x_\mu)^{d_s}$ ($k=1,2,...,i$). Since
$J_1$ is Borel fixed, assume by induction that for some $1\leq
k<i$, the ideal $J_1+...+J_k$ is Borel fixed and let $\mathbf{n}
\in J_1+...+J_k+J_{k+1}$. If $\mathbf{n} \in J_1+...+J_k$, we are
done, so assume that $\mathbf{n} \in G(N_{k+1})
\setminus(J_1+...+J_{k})$. Then write
\[
    \mathbf{n}=\mathbf{n'}\mathbf{n''}
    \]
for some $\mathbf{n'} \in G(N_{k+1})$ and $\mathbf{n''} \in
G\left((x_{k+1},...,x_{\mu})^{d_{s}}\right)$.  Now observe that
$x_{k+1}$ must divide $\mathbf{n'}$ because $\mathbf{n} \notin
J_1+...+J_{k}$. Next, if $r<t$ and $x_t$ divides $\mathbf{n}$,
then we see that  $\mathbf{n}_{t \to r}$ is in $J_{k+1}$. Indeed,
it is easy to verify this when $x_t$ divides $\mathbf{n'}$,
because $N_{k+1}$ is Borel fixed, so assume that  $x_t$ does not
divide $\mathbf{n'}$. Then $x_t$ divides $\mathbf{n''}$,
  so $t\geq k+1$. If $k+1\leq r$, then we have
  \[
   \frac{\mathbf{n}x_{r}}{x_{t}}= \mathbf{n'} \frac{\mathbf{n''}x_{r}}{x_t}\in N_{k+1}(x_{k+1},...,x_{\mu})^{d_{s}}=J_{k+1},
  \]
  while if $r<k+1$,
\[
   \frac{\mathbf{n}x_{r}}{x_{t}}=\frac{\mathbf{n'}x_{r}}{x_{k+1}}\cdot \frac{\mathbf{n''}x_{k+1}}{x_t}\in
   N_{k+1}(x_{k+1},...,x_{\mu})^{d_{s}}=J_{k+1},
  \]
  because $N_{k+1}$ is Borel fixed. Thus
 \[
 \dfrac{\mathbf{n}x_{r}}{x_{t}} \in J_1+...+J_k+J_{k+1}
 \]
  in all cases and the proof of part (a) is complete. For part
  (b), let $\mathbf{m} \in G(N_j(x_1,...,x_\mu)^{d_s})$
  and $\mathbf{n} \in G(N_{j+1}(x_{j+1},...,x_\mu)^{d_s})$, and
  write $\mathbf{m}=\mathbf{m}_1\mathbf{m}_2$ and
  $\mathbf{n}=\mathbf{n}_1\mathbf{n}_2$, where $\mathbf{m}_1\in
  G(N_j)$, $\mathbf{m}_2 \in G((x_1,...,x_\mu)^{d_s})$, $\mathbf{n}_1\in
  G(N_{j+1})$ and $\mathbf{n}_2 \in G((x_{j+1},...,x_\mu)^{d_s})$.
  Then, note that $\mathbf{m}_1\mathbf{n}_2$ divides $lcm(\mathbf{m},\mathbf{n})$. This implies that
  $lcm(\mathbf{m},\mathbf{n})$ is in
  $N_{j}(x_{j+1},...,x_\mu)^{d_s}$, and so
  \[
    N_j(x_1,...,x_\mu)^{d_s} \cap   N_{j+1}(x_{j+1},...,x_\mu)^{d_s}\subseteq
    N_{j}(x_{j+1},...,x_\mu)^{d_s}.
  \]
The opposite containment is obvious, so the proof of part (b) is
complete.
\endproof
\begin{rem}
 Part (a) with $i=\lambda_{s-1}$ and $\mu=\lambda_s$ yields
  \[
  I=N_1
  (x_1,...,x_{\lambda_s})^{d_s}+N_2(x_2,...,x_{\lambda_s})^{d_s}+...+N_{\lambda_{s-1}}(x_{\lambda_{s-1}},...,x_{\lambda_s})^{d_s}.
\]
\end{rem}

\begin{exs} 1) For the ideal $I=<bd^2>$ in $\Bbbk[a,b,c,d]$, we have
$s=2$, $\lambda_1=2$, $d_1=1$, $\lambda_2=4$ and $d_2=2$.
Moreover, $N_1=<a>$, $N_2=<b>=(a,b)$. Therefore,
\[
 I=<a>(a,b,c,d)^2+<b>(b,c,d)^2
 \]
 2) For the ideal $I=<bcd>$ in $\Bbbk[a,b,c,d]$, we have
$s=3$, $\lambda_1=2$, $\lambda_2=3$, $\lambda_3=4$ and
$d_1=d_2=d_3=1$. Moreover, $N_1=<a^2>$, $N_2=<b^2>$ and
$N_3=<bc>$. Therefore,
\[
 I=<a^2>(a,b,c,d)+<b^2>(b,c,d)+<bc>(c,d).
 \]
 \end{exs}

\begin{figure}[htbp]
\begin{center}
\input{example51.pstex_t}
\end{center}
\end{figure}
\
\\
\\
 \textbf{Proof of Theorem \ref{P2:CR}}: By induction on $s$. For $s=1$, we
are done. Assume that $s>1$ and that we have obtained $Q\left(\dis
\prod_{j=1}^{k} I_{\lambda_j}^{d_j}\right)$ for all $k<s$. By the
inductive hypothesis, there is a polyhedral cell complex $Q(N_i)$
that supports a minimal free resolution for the ideals $N_{i}$,
for all $1 \leq i\leq \lambda_{s-1}$. Moreover, from our
construction we see that $Q(N_i)$ is a subcomplex of $Q(N_{i+1})$
for all $1\leq i<\lambda_{s-1}$. Set
$J_i=N_k(x_i,...,x_{\lambda_s})^{d_s}$
($i=1,2,...,\lambda_{s-1}-1$). From lemmas \ref{L2:CR},
\ref{L3:CR} and \ref{L4:CR}  it follows that the polyhedral cell
complexes $C_i$ and $D_i$ ($i=1,2,...,\lambda_s$) defined by
\[
   C_i:=Q(N_i)\times P_{d_{s}}(x_i,...,x_{\lambda_s})
\]
and
\[
 D_i:=C_i \cap C_{i+1}= Q(N_i)\times P_{d_{s}}(x_{i+1},...,x_{\lambda_s}).
\]
 support a  minimal free resolution of $J_i$ and $(J_1+...+J_i) \cap J_{i+1}$,
 respectively, for all $1\leq i<\lambda_{s-1}$. Thus, lemma \ref{L1:CR} implies
 that the polyhedral cell complex $C'_k$, which is defined recursively by
 \[
  C'_1=C_1,\quad \text{and} \quad C'_{i+1}=C'_{i} \cup C_{i+1}
 \]
 for $1\leq i< \lambda_{s-1}$, supports a (minimal) free resolution of
 $J_1+J_2+...+J_i$.
 Accordingly, set
\[
  Q(I):=C'_{\lambda_{s-1}}= C_1 \cup C_2 \cup \cdots \cup C_{\lambda_{s-1}}
\]
and the construction of our polyhedral cell complex is done by
induction. Also, from our construction it follows that
$Q(\mathbf{m})$ is a subcomplex of $P_d(x_1,...,x_n)$, where
$d=degree(\mathbf{m})$, as desired. Finally, it is easy to see as
in \ref{P1:CR} that $Q(\mathbf{m})$ is also shellable in order to
complete the proof.\endproof

\subsection{$d$-generated Borel fixed ideals}

 Next we would like to generalize theorem \ref{P2:CR} to
the case of any Borel fixed ideal generated in one degree. Before
we prove this in \ref{P3:CR}, we need more preliminary results.
Recall that for two monomials $\mathbf{m}_1$ and $\mathbf{m}_2$ of
the same degree, $\mathbf{m}_1 \succ_{rlex} \mathbf{m}_2$ means
that the rightmost non-zero entry of the difference
$e(\mathbf{m}_1)-e(\mathbf{m}_2)$ is negative.

\begin{lemma} \label{L5:CR}Let $\mathbf{m}_1$ and $\mathbf{m}_2$ be
two monomials of the same degree $d$, such that $\mathbf{m}_1
\succ_{rlex} \mathbf{m}_2$, which minimally generate in the Borel
sense an ideal
\[
  I =<\mathbf{m}_1, \mathbf{m}_2>.
\]
Then
\[
   <\mathbf{m}_1> \cap <\mathbf{m}_2>
\]
is a principal Borel ideal. Moreover,
\[
   Q(\mathbf{m}_1) \cap Q(\mathbf{m}_2)=Q(\mathbf{m}).
\]
\end{lemma}
    \textbf{Proof}. First assume that $\mathbf{m}_1=x_1^{a_1}x_2^{a_2}\cdots x_n^{a_n}$ and
$\mathbf{m}_2=x_1^{b_1}x_2^{b_2}\cdots x_n^{b_n}$. Then define
$\mu_n=min\{a_n,b_n\}$ and the natural numbers $\mu_i$ for
$i=n-1,...,1$ recursively, by setting
\[
   \mu_{i}=min\{a_{i}+...+a_n, b_{i}+...+b_n\}-(\mu_{i+1}+...+\mu_n).
   \]
 Define the following monomial of degree $d$
\[
  \mathbf{MIN}(\mathbf{m}_1,\mathbf{m}_2):=x_1^{\mu_1}x_2^{\mu_2}\cdots x_n^{\mu_n}
\]
 From our choice of the $\mu_i$'s, we have
\[
  \mu_{n-i}+ \mu_{n-i+1}+...+\mu_n\leq a_{n-i}+a_{n-i+1}+...+a_n,
\]
and
\[
  \mu_{n-i}+\mu_{n-i+1}+...+\mu_n  \leq b_{n-i}+b_{n-i+1}+...+b_n ,
\]
for all $i=0,1,...,n-1$. Therefore,
\[
  <\mathbf{MIN}(\mathbf{m}_1,\mathbf{m}_2)> \subseteq <\mathbf{m}_{1}> \cap <\mathbf{m}_2>.
\]
Now let $\mathbf{m}=x_1^{c_1}x_2^{c_2}\cdots x_n^{c_n}$ be in
$G(<\mathbf{m}_{1}>)$ and  let
$\mathbf{n}=x_1^{d_1}x_2^{d_2}\cdots x_n^{d_n}$ be in
$G(<\mathbf{m}_2>)$. We want to show that
$lcm(\mathbf{m},\mathbf{n})$ is in
$<\mathbf{MIN}(\mathbf{m}_1,\mathbf{m}_2)>$, so we may assume that
$c_1<d$ and $d_1<d$. Next, let $k$ be the largest positive integer
such that
\[
max\{c_1,d_1\}+...+max\{c_k,d_k\}<d.
\]
Then set $\nu_i=max\{c_i,d_i\}$ for $i=1,2,...,k$ and
$\nu_{k+1}=d-(\nu_1+...+\nu_k)$. Since
\begin{align}\notag
 \nu_1+...+\nu_i &\geq max\{c_1+...+c_i,d_1+...+d_i\}\\\notag
             &\geq max\{a_1+...+a_i,b_1+...+b_i\},\notag
\end{align}
for all $1\leq i\leq k$, we see that
\begin{align}\notag
  \nu_{i+1}+...+\nu_{k+1} &\leq  d-max\{a_1+...+a_i,b_1+...+b_i\}\\\notag
                    &=min\{a_{i+1}+...+a_n,b_{i+1}+...+b_{n}\}\\\notag
                    &=\mu_{i+1}+...+\mu_n.\notag
\end{align}
Therefore, the monomial $x_1^{\nu_1}x_2^{\nu_2}\cdots
x_{k+1}^{\nu_{k+1}}$ is a minimal generator of
$<\mathbf{MIN}(\mathbf{m}_1,\mathbf{m}_2)>$ and divides
$lcm(\mathbf{m},\mathbf{n})=x_1^{max\{c_1,d_1\}}x_2^{max\{c_2,d_2\}}\cdots
x_{n}^{max\{c_n,d_n\}}$. Therefore, $lcm(\mathbf{m},\mathbf{n})$
is in $<\mathbf{MIN}(\mathbf{m}_1,\mathbf{m}_2)>$, and so
\[
   <\mathbf{m}_{1}> \cap <\mathbf{m}_2> \subseteq \,<\mathbf{MIN}(\mathbf{m}_1,\mathbf{m}_2)>.
\]
Thus, the proof of our first claim is complete with
$\mathbf{m}:=\mathbf{MIN}(\mathbf{m}_1,\mathbf{m}_2)$. Now note
that $Q(\mathbf{m}_{1}) \cap Q(\mathbf{m}_2)$ is the union of all
the convex polytopes of the polyhedral cell complex
$P_d(x_1,...,x_n)$ with vertices in $<\mathbf{m}_{1}> \cap
<\mathbf{m}_2>=<\mathbf{m}>$. Since $Q(\mathbf{m})$ is the union
of all the convex polytopes of the polyhedral cell complex
$P_d(x_1,...,x_n)$ with vertices in $<\mathbf{m}>$, we must have
\[
   Q(\mathbf{m}_1) \cap Q(\mathbf{m}_2)=Q(\mathbf{m}).
\]
\endproof

\begin{exs}\label{Exs:CR} 1) Let $\mathbf{m}_1=b^5c$, $\mathbf{m}_2= ab^3c^2$ and
$\mathbf{m}_3= a^2c^4$ in $\Bbbk[a,b,c]$. Then
\begin{align}\notag
   <\mathbf{m}_1>\cap <\mathbf{m}_2> &= <b^5c> \cap <ab^3c^2>=<ab^4c>,\\\notag
   <\mathbf{m}_1>\cap <\mathbf{m}_3> &= <b^5c> \cap <a^2c^4>=<a^2b^3c>\notag
   \end{align}
Also,
\begin{align}\notag
   <\mathbf{m}_1>\cap <\mathbf{m}_2,\mathbf{m}_3> &=<ab^4c,a^2b^3c>=<ab^4c>.\notag
\end{align}

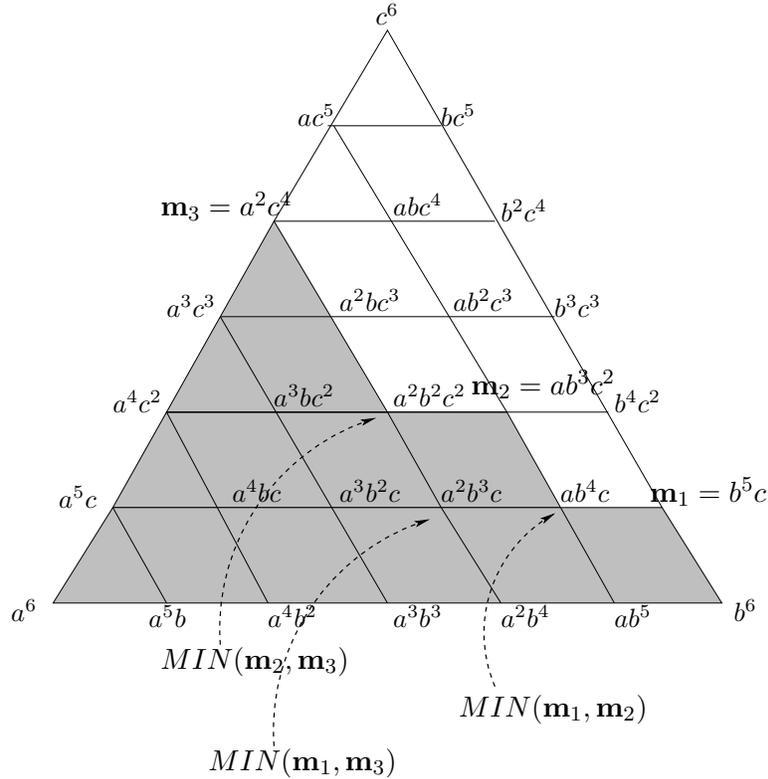
\begin{figure}[htbp]
\begin{center}
\input{example71.pstex_t}
\end{center}
\end{figure}

In general, the intersection of a principal Borel ideal with a
non-principal Borel ideal is not principal. For example, in
$\Bbbk[a,b,c,d]$ we have
\[
   <ab^4c^3d>\cap <a^2b^4cd^2,a^3bc^2d^3> =<a^3b^2c^3d, a^2b^4c^2d>.
\]
\end{exs}

The following lemma is sufficient for our purposes.

\begin{lemma} \label{L6:CR} Let $\mathbf{m}_1$, $\mathbf{m}_2$, ...,$\mathbf{m}_s$  be $s$ monomials of the same degree
$d$, such that $\mathbf{m}_1 \succ_{rlex} \mathbf{m}_2 \succ_{rlex}
...\succ_{rlex} \mathbf{m}_s$, which minimally generate in the Borel
sense an ideal
\[
    I=<\mathbf{m}_1,...,\mathbf{m}_s>.
    \]
     Then for all $j<s$,
\[
   <\mathbf{m}_j> \cap <\mathbf{m}_{j+1},...,\mathbf{m}_s>
\]
is a Borel fixed ideal, which is minimally generated in the Borel
sense by \emph{at most} $s-j$ monomials
$\mathbf{n}_{j+1},...,\mathbf{n}_{s}$ of degree $d$, with
$\mathbf{n}_k \succ_{rlex} \mathbf{m}_j$ for $k=j+1,...,s$.
\end{lemma}
\textbf{Proof}. First note that the case where $j=s-1$ is
essentially the previous lemma  \ref{L5:CR}. Now for all $j<s$, we
have
\begin{align}\notag
    <\mathbf{m}_{j}> \cap <\mathbf{m}_{j+1},...,\mathbf{m}_s>&=\left(<\mathbf{m}_j> \cap <\mathbf{m}_{j+1}>\right)+...+\left(<\mathbf{m}_j> \cap <\mathbf{m}_{s}>\right)\\\notag
                                        &=<\mathbf{MIN}(\mathbf{m}_{j},\mathbf{m}_{j+1})> +...+ <\mathbf{MIN}(\mathbf{m}_{j},\mathbf{m}_{s})>\\\notag
                                        &=<\mathbf{n}_{j+1},...,\mathbf{n}_s>\notag
\end{align}
where
\[
  \mathbf{n}_k:=\mathbf{MIN}(\mathbf{m}_j,\mathbf{m}_{k})
\]
for $k=j+1,...,s$. As we saw in example \ref{Exs:CR}, some of the
$\mathbf{n}_k$'s might be redundant, so the above intersection is
minimally generated in the Borel sense by \emph{at most} $s-j$
monomials $\mathbf{n}_{j+1},...,\mathbf{n}_{s}$ of degree $d$,
with $\mathbf{n}_k
 \succ_{rlex} \mathbf{m}_j$ for $k=j+1,...,s$.\endproof

Now we are ready to prove our most general result.

\begin{thm} \label{P3:CR} Let $\mathbf{m}_1$, $\mathbf{m}_2$, ...,$\mathbf{m}_s$  be $s$ monomials of the same degree
$d$, such that $\mathbf{m}_1 \succ_{rlex} \mathbf{m}_2 \succ_{rlex}
...\succ_{rlex} \mathbf{m}_s$, which minimally generate in the Borel
sense an ideal
\[
    I=<\mathbf{m}_1,...,\mathbf{m}_s>.
    \]
Then there exists a polyhedral cell complex
$Q(\mathbf{m}_1,...,\mathbf{m}_s)$ that supports a minimal free
resolution of $I$. Moreover, $Q(\mathbf{m}_1,...,\mathbf{m}_s)$ is
the union of all the convex polytopes of the polyhedral cell complex
$P_d(x_1,...,x_n)$ with vertices in
$I=<\mathbf{m}_1,...,\mathbf{m}_s>$.
\end{thm}
\textbf{Proof}: For $s=2$ both of our claims were proved in Lemma
\ref{L5:CR}.  So assume that $s> 2$, and for all $j<s$ set
\[
  I_j=<\mathbf{m}_j,...,\mathbf{m}_s>
\]
 Next suppose that for some $j< s$ we have
constructed a polyhedral cell complex $Q(K)$ that supports a
minimal free resolution of any Borel fixed ideal $K$, which is
minimally generated in the Borel sense by at most $s-j$ monomials
of the same degree $d$. Assume also that $Q(K)$ is the union of
all the convex polytopes of the polyhedral cell complex
$P_d(x_1,...,x_n)$ with vertices in $K$.

From lemma  \ref{L6:CR}, we see that
\[
   <\mathbf{m}_j> \cap <\mathbf{m}_{j+1},...,\mathbf{m}_s>=<\mathbf{n}_{j+1},...,\mathbf{n}_{s}>
\]
is a Borel fixed ideal, which is minimally generated by \emph{at
most} $s-j$ monomials $\mathbf{n}_{j+1},...,\mathbf{n}_{s}$ of
degree $d$. Thus, so far we have constructed the polyhedral cell
complex $Q(\mathbf{m}_j)$ in theorem \ref{P2:CR}, and the
polyhedral cell complexes $Q(\mathbf{m}_{j+1},...,\mathbf{m}_s)$
and $Q(\mathbf{n}_{j+1},...,\mathbf{n}_{s})$, by the inductive
hypothesis. Moreover, by the inductive hypothesis,
$Q(\mathbf{m}_j)\cap Q(\mathbf{m}_{j+1},...,\mathbf{m}_s)$ is the
union of all the convex polytopes of the polyhedral cell complex
$P_d(x_1,...,x_n)$ with vertices in $<\mathbf{m}_j>\cap
<\mathbf{m}_{j+1},...,\mathbf{m}_s>$. Since
$Q(\mathbf{n}_{j+1},...,\mathbf{n}_{s})$ is the union of all the
convex polytopes of the polyhedral cell complex $P_d(x_1,...,x_n)$
with vertices in $<\mathbf{n}_{j+1},...,\mathbf{n}_{s}>$, we must
have
\[
   Q(\mathbf{m}_j)\cap Q(\mathbf{m}_{j+1},...,\mathbf{m}_s)=Q(\mathbf{n}_{j+1},...,\mathbf{n}_{s}).
\]
Since the rest of the hypotheses of lemma  \ref{L1:CR} are easily
checked to be satisfied, we conclude that the complex
\[
  X_j:=Q(\mathbf{m}_j) \cup Q(I_{j+1})=Q(\mathbf{m}_j)\cup Q(\mathbf{m}_{j+1})\cup ... \cup Q(\mathbf{m}_s)
\]
supports a minimal free resolution of the ideal $I_j$. Thus
\[
  X:=X_1=Q(\mathbf{m}_1)\cup Q(\mathbf{m}_{2})\cup ... \cup Q(\mathbf{m}_s).
\]
supports a minimal free resolution of $I_1=I$.\endproof

\section{The lcm-lattice } \label{S:LCM}

 The lcm-lattice of an arbitrary monomial ideal $I$ was introduced
 in \cite{GPW}, where the authors show how its structure relates to the
 Betti numbers and the maps in the minimal free resolution of $I$. The lcm-lattice of $I$, with $G(I)=\{\mathbf{m}_1,\mathbf{m}_2,...,\mathbf{m}_r\}$,
 is denoted by $L_I$. This is the lattice with elements labeled by the least common
multiple of $\mathbf{m}_1,\mathbf{m}_2,...,\mathbf{m}_r$ ordered
by divisibility; that is, if $\mathbf{n}$ and $\mathbf{m}$ are
distinct elements of $L_I$, then $\mathbf{m}\prec \mathbf{n}$ if
and only if $\mathbf{m}$ divides $\mathbf{n}$. Moreover, we
include $\hat{0}:=1$ as the bottom element, while
$\hat{1}=lcm(\mathbf{m}_1,\mathbf{m}_2,...,\mathbf{m}_r)$ is the
top element. We say that $\mathbf{n}$ \emph{covers} $\mathbf{m}$
and we write $\mathbf{m} \to \mathbf{n}$, if $\mathbf{m}\prec
\mathbf{n}$ and if there is no element $\mathbf{k}\ne \mathbf{n},
\mathbf{m}$ of $L_I$ such that $\mathbf{m}\prec \mathbf{k} \prec
\mathbf{n}$.

We would like to find a labelling of the edges of $L_I$ with the
following property: for all elements $\mathbf{m}$ and $\mathbf{n}$
in $L_I$ with $\mathbf{m}\prec \mathbf{n}$, there exists a unique
\emph{increasing} maximal chain from  $\mathbf{m}$ to $\mathbf{n}$
and it is lexicographically strictly first than all other maximal
chains from  $\mathbf{m}$ to $\mathbf{n}$. This would prove that
$L_I$ is shellable (see \cite{BjornerWachs1},
\cite{BjornerWachs2}) in a way different than \cite{Ekki}. Finding
such a labelling is still an open problem.

\begin{rem}
The natural labelling which assigns to each edge $\mathbf{m} \to
\mathbf{n}$ the integer
$\text{max}\left(\frac{\mathbf{n}}{\mathbf{m}}\right):=\text{max}\{i|x_i\quad
\text{divides}\quad
    \frac{\mathbf{n}}{\mathbf{m}}\}$ does not
work.
\begin{ex}
Let
\begin{align}\notag
    I&=<ab,ac,ad^2,b^2cd^2>\\\notag
     &=(a^2, ab, b^5 , ac, b^4c, b^3c^2 , b^2c^3, b^4d, b^3cd, b^2c^2d, ad^2 , b^3d^2, b^2cd^2
     )\notag
\end{align}
The interval $[1,ab^2cd^2]$ of $L_I$ is
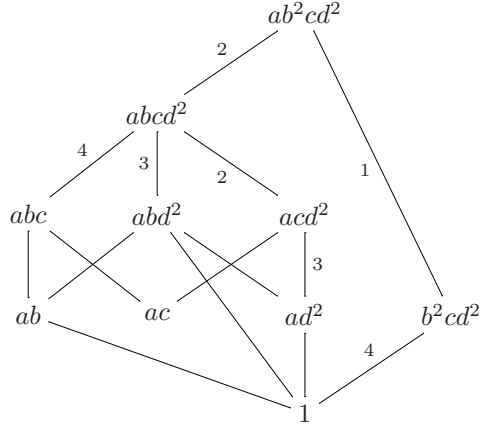
\[
 \xymatrix{
           &                & {ab^2cd^2 }      &             & &\\
           &       {abcd^2}   \ar@{-}[ur]^2        &       &             & &\\
             {abc} \ar@{-}[ur]^4        & {abd^2}\ar@{-}[u]^3&  {acd^2}\ar@{-}[ul]^2      &            & &\\
             {ab}\ar@{-}[u]\ar@{-}[ur]  & {ac}\ar@{-}[ul] \ar@{-}[ur]   & {ad^2 }\ar@{-}[ul]\ar@{-}[u]_3     &  {b^2cd^2}\ar@{-}[uuul]^1   &  &\\
           &            & {1} \ar@{-}[ull] \ar@{-}[uul]\ar@{-}[u]
           \ar@{-}[ur]^4
                                                                       &              &&}
\]
Hence there is no decreasing sequence of labels from $ab^2cd^2$ to
1 (or even to $abc$).
\end{ex}
\end{rem}

 The above example shows also that the lcm lattice of a
Borel fixed ideal need not be ranked in general. However, if $I$
is generated in the same degree then we prove the following.

\begin{prop}\label{K:LCM} The lcm-lattice $L_I$ of a $d$-generated Borel
fixed ideal $I$ is ranked.
 \end{prop}
\textbf{Proof}: Let
$I=(\mathbf{m}_1,\mathbf{m}_2,...,\mathbf{m}_r)$ be minimally
generated by $\mathbf{m_1},\mathbf{m_2},...,\mathbf{m_r}$ in the
same degree $d$ and let $\mathbf{m} \ne \hat{1}=
lcm(\mathbf{m}_1,\mathbf{m}_2,...,\mathbf{m}_r)$ be in the
lcm-lattice $L_I$ of $I$. Assume that
$\mathbf{m}=lcm(\mathbf{m}_{\alpha},\mathbf{m}_{\beta},...,\mathbf{m}_{\gamma})$,
with
\[
    e(\mathbf{m}_{\alpha})=(a_1,a_2,...,a_n),\quad e(\mathbf{m}_{\beta})=(b_1,b_2,...,b_n) \quad \text{...} \qquad e(\mathbf{m}_{\gamma})=(c_1,c_2,...,c_n).
\]
In order to show that the lattice is ranked it suffices to prove
that $deg(\mathbf{n})=1+deg(\mathbf{m})$ for all $\mathbf{n}$ that
cover $\mathbf{m}$. There exists a $\mathbf{m}_{\delta}$ in $I$,
with $e(\mathbf{m}_{\delta})=(d_1,d_2,...,d_n)$ such that
$\mathbf{n}=lcm(\mathbf{m}_{\alpha},\mathbf{m}_{\beta},...,\mathbf{m}_{\gamma},
\mathbf{m}_{\delta})=lcm(\mathbf{m},\mathbf{m}_{\delta})$. Also,
there is at least one $j$ such that $d_j>max\{a_j,b_j,...,c_j\}$.
Without loss of generality assume that for that $j$, it is
$max\{a_j,b_j,...,c_j\}=a_j$. If there is some $k$ with $j<k\leq
n$ such that $a_k \ne 0$, then
\[
    \ell:=lcm ((\mathbf{m}_{\alpha})_{k \to j},\mathbf{m}_{\alpha},\mathbf{m}_{\beta},...,\mathbf{m}_{\gamma})
\]
has degree $deg(\mathbf{\ell})=1+deg(\mathbf{m})$, divides
$\mathbf{n}$ and is divisible by $\mathbf{m}$. The minimality of
$\mathbf{n}$ forces $\mathbf{\ell}=\mathbf{n}$ and so
$deg(\mathbf{n})=1+deg(\mathbf{m})$. Now assume that $a_k=0$ for
all $j<k\leq n$. Then, there is an $i<j$ such that
$max\{a_i,b_i,...,c_i\}>d_i$. [Indeed, suppose to the contrary
that $d_i \geq max\{a_i,b_i,...,c_i\}$ for all $i<j$. Then,
\[
  d\geq \sum_{i=1}^j d_i > \sum_{i=1}^j max\{a_i,b_i,...,c_i\}
  \geq \sum_{i=1}^j a_i = \sum_{i=1}^n a_i =d,
\]
a contradiction.] Then
\[
    \ell:=lcm ((\mathbf{m}_{\delta})_{j \to i},\mathbf{m}_{\alpha},\mathbf{m}_{\beta},...,\mathbf{m}_{\gamma})
\]
has degree $deg(\mathbf{n})-1$, divides $\mathbf{n}$ and is
divisible by $\mathbf{m}$. Hence, $\mathbf{\ell}=\mathbf{m}$ and so
$deg(\mathbf{n})=1+deg(\mathbf{m})$, as desired. The proof is
complete.\endproof

\begin{rems} 1) The above proof applies with minor modifications to the
case of a \emph{strongly stable square-free} ideal generated in
the same degree. A monomial ideal $I$ is called \emph{strongly
stable square-free} if all monomials in $G(I)$ are square-free and
for every $\mathbf{m}$ in $G(I)$, if $x_t$ divides $m$ and $x_s$
does not divide $\mathbf{m}$ ($1\leq s<t$), then $\mathbf{m}_{t
\to s}$ is in $I$.

2) There exists a $d$-generated Borel fixed ideal
$I=(\mathbf{m}_1,\mathbf{m}_2,...,\mathbf{m}_r)$ minimally
generated by $\mathbf{m}_1,\mathbf{m}_2,...,\mathbf{m}_r$ and an
element $\mathbf{m}$ of $L_I$ of degree $d+1$, such that for some
$1\leq s <t\leq n$,
\begin{itemize}
    \item[(i)] $x_t$ divides $\mathbf{m}$
    \item[(ii)]$x_s^{d_s}$ does not divide $\mathbf{m}$

    \end{itemize}
where $d_i$ is the largest positive integer such that
$x_{i}^{d_i}$ divides
$lcm(\mathbf{m}_1,\mathbf{m}_2,...,\mathbf{m}_r)$, and
\begin{itemize}
    \item[(iii)] $\mathbf{m}_{t \to s}$ is not in $L_I$,
\end{itemize}

 \begin{ex} Let
\[
  I=<x_1 x_3^3, x_2^2 x_3 x_4>.
\]
 Then $d_3=3$ and $x_2^2 x_3^2 x_4=lcm(x_2^2 x_3 x_4,x_2^2
 x_3^2)$ is in $L_I$, but $x_2^2 x_3^3=(x_2^2 x_3^2 x_4)_{4\to 3}$
 is \emph{not} in $L_I$.
 \end{ex}
\end{rems}

\textbf{Acknowledgements}: I would like to express my thanks to
Mike Stillman for his support and for the numerous discussions we
have had on this project. I would like to thank Irena Peeva for
her suggestions on section 4 and the discussions we had on them.
Also, I would like to thank Volkmar Welker and Ekki Batzies for
valuable e-mail communications, as well as Victor Reiner and
Mauricio Velasco for helpful discussions.

Address:

 Cornell University, Department of Mathematics,
 Ithaca, NY 14853,

 asin@math.cornell.edu

\end{document}